\documentclass{article}
\usepackage{amssymb,amsmath}
 
\input xy 
\xyoption{all}
 
\begin{document}
\title{Hodge theory for elliptic complexes over unital  $C^*$-algebras}
\author{Svatopluk Kr\'ysl \footnote{{\it E-mail address:} Svatopluk.Krysl@mff.cuni.cz, {\it Tel./Fax:} $+$ 420  222 323 221/ $+$ 420 222 323 394} \\{\it \small 
Charles University in Prague, Faculty of Mathematics and Physics
}}
\maketitle \noindent
\centerline{\large\bf Abstract}

 For a unital Banach $C^*$-algebra $A,$  we prove that the cohomology groups of  $A$-elliptic complexes of pseudodifferential  operators in finitely generated projective $A$-Hilbert 
bundles over  compact manifolds are norm complete topological vector spaces and finitely generated $A$-modules provided the images of certain extensions of the so called associated Laplacians are closed. 
This establishes a Hodge type theory for these structures.


{\it Math. Subj. Class.:} 35J05, 46L87,  58A14, 58J10
  
\section{Introduction}

  In this paper, we  deal  with fields having their values in possibly infinite rank vector bundles and systems of pseudodifferential equations for these fields. 
To be more precise, we focus our attention at complexes of pseudodifferential operators acting between smooth sections of finitely generated projective $A$-Hilbert bundles over compact manifolds.
We recall  a definition  of  an $A$-elliptic complex and prove a Hodge type theory for a certain subclass of  them. By a Hodge type theory, we mean not only constructing an $A$-linear isomorphism between the
cohomology groups and the spaces of harmonic elements of the complex, but also  giving a description of the  the cohomology groups from a topological  point of view.

Let $A$ be a unital  $C^*$-algebra. 
A finitely generated projective $A$-Hilbert bundle is, roughly speaking, a fiber bundle the total space of which is a Banach  manifold, and the fibers of which  are finitely 
generated projective Hilbert 
$A$-modules. 
A Hilbert $A$-module is a module over $A$
which is first, equipped with a positive definite $A$-sesquilinear 
map   with values in $A,$   and second, it is a  complete topological vector space with respect to the norm derived from the $A$-sesquilinear map and  from the Banach 
norm in $A.$ 
Our reference for general Hilbert $A$-modules is Lance \cite{Lance}, and 
for the finitely generated and projective ones  the monograph of Solovyov, Troitsky \cite{TM}.

One of the basic basic steps in proving the so called Mishchenko-Fomenko index theorem (Fomenko, Mishchenko \cite{FM}) for 
$A$-elliptic pseudodifferential operators acting on sections 
of finitely generated projective $A$-Hilbert bundles over a given compact manifold  is a construction of  certain pseudoinverses
to extensions of such operators to  the Sobolev type completions of the space of smooth sections. See Fomenko, Mishchenko  \cite{FM} for this construction.
In \cite{FM}, not only  a paramatrix for $A$-elliptic pseudodifferential operators is constructed, but the authors also prove that such operators 
are $A$-Fredholm. In particular, the kernels of such operators are finitely generated projective Hilbert $A$-modules. Let us notice,  that we use the   theorems and notions
 mentioned in this paragraph  in the form in which they appear in Solovyov, Troitsky \cite{TM}.

For a chain complex $(D_i, \Gamma^i)_{i\in \mathbb{N}_0}$   of pre-Hilbert $A$-modules and adjointable pre-Hilbert $A$-module homomorphisms,
we can form a sequence of the so called {\it associated Laplacians} $\triangle_i,$ $i\in \mathbb{N}_0,$ without supposing any topology on the modules. 
Namely, one sets $\triangle_i = D_i^* D_i +D_{i-1} D_{i-1}^*$ simply.
Assume that each $\triangle_i$ possesses a parametrix, i.e., there exist pre-Hilbert $A$-module homomorphisms $g_i: \Gamma^i \to \Gamma^i$ and
 $p_i: \Gamma^i \to \Gamma^i$ such that
parametrix equations $1 = \triangle_ig_i + p_i,$  $1 = g_i \triangle_i + p_i$ hold, and  $p_i$ maps $\Gamma^i$ into the kernel of  $\triangle_i.$
We call complexes with such behaved Laplacians {\it parametrix possessing}.
It is quite interesting, and we prove that  in this case, the pseudoinverses   $g_i$ are necessarily chain homomorphisms, i.e.,
$g_{i+1} D_i = D_i g_i.$ This is at least implicitly known also in the
finite rank situation, but we derive this property in a purely  algebraic way. 
 Knowing this, it is not hard to show that the cohomology 
groups of a parametrix possessing complex are $A$-linearly isomorphic to the kernels of the appropriate Laplacians, i.e., to the spaces of {\it harmonic elements.} This establishes an abstract 
or, say, an algebraic Hodge theory. At this level, a topological characterization of the cohomologies  is missing.
 
 By an {\it $A$-elliptic  complex}, we mean a complex of  pseudodifferential operators acting on smooth sections of  $A$-Hilbert bundles 
  the associated symbol sequence of which is exact
out of the zero section of the cotangent bundle.
We prove  that the symbols of the  Laplacians associated to an $A$-elliptic   complex  
 are isomorphisms out of the zero section. From now on till the end of this paragraph, 
all mentioned bundles are supposed to be finitely generated projective $A$-Hilbert bundles over compact 
manifolds.
We use a generalization of the Sobolev embedding theorem for sections of $A$-Hilbert 
bundles, which we prove in this text, to derive a
regularity of  $A$-elliptic pseudodifferential operators.  
  If we moreover, assume that the images of  certain extensions, specified in the paper, of each of the associated Laplacians of an $A$-elliptic complex 
 are closed,  we are able to construct parametrix equations for the Laplacians using the regularity.
Thus,  we prove that such $A$-elliptic complexes   are parametrix possessing.
Especially, their cohomology groups are isomorphic to the spaces of harmonic elements   as $A$-modules. This establishes the algebraic Hodge theory for these 
 structures. Consequently, the mentioned  $A$-Fredholm property of $A$-elliptic pseudodifferential operators implies  that the cohomology groups of this subclass of $A$-elliptic  complexes are
 finitely generated  $A$-modules. Using some basic topological reasoning, we prove that the cohomologies of these complexes are Banach topological vector spaces.

The regularity of an $A$-elliptic
operator in finitely generated projective $A$-Hilbert bundles over compact manifolds was already proved, e.g.,  in \cite{TM} (Theorem 2.1.145).
We don't give a substantially new proof of this assertion, but we try to write an as much self-contained proof of this fact as possible. 
We also  notice  that although we could have considered general  $C^*$-algebras at least until Lemma 5 inclusively, we  decided to suppose  that all  $C^*$-algebras 
are unital from the beginning.

 Further, let us
remark  that there exist generalizations of the classical Hodge theory in directions different from that one described here.
See, e.g., Bartholdi et al. \cite{B} and Smale et al. \cite{Smale} for a generalization to complete separable metric spaces endowed with a probabilistic measure, 
and also for further references given there. We do not attempt to give a full reference to the  topic of complexes of pseudodifferential operators in $A$-Hilbert bundles, and refer the reader to 
Troitsky \cite{Troi}, Solovyov, Troitsky \cite{TM} and Schick \cite{Schick}.
We develop the presented theory mainly in order to enable a description of solutions to equations for operator or Hilbert module 
valued fields similar to the equations  appearing in  Quantum field theory. Our motivation comes, however, from geometric and deformation quantization via 
the so called {\it symplectic spinor fields}. See Kostant \cite{Kostant}, Fedosov \cite{Fedosov} and Habermann \cite{KH} for this context. 

   In the second section,  we prove a theorem on the homotopy properties  of the  parametrix possessing complexes (Theorem 3), and in Theorem 4, 
we derive the algebraic form of the Hodge theory for them. In the third part,  we recall the definition of a finitely generated projective  $A$-Hilbert bundle, Sobolev type completions of smooth sections,
 and  the definition of the Fourier transform in $A$-Hilbert bundles. The embedding theorem is stated as Lemma 5. The regularity for  $A$-elliptic 
operators is proved in Theorem 7,  and the smooth pseudoinverses for the $A$-elliptic and self-adjoint ones are constructed in Theorem 8.
 In the fourth section, we recall the notion of an $A$-elliptic complex and prove that its cohomology groups are finitely generated   $A$-modules and
Banach topological vector spaces   under the   mentioned condition on the images of the extensions of the associated Laplacians  (Theorem 11).
 

{\bf Preamble:} In the whole text  when not said otherwise, manifolds, fibrations  (bundle projections, total and base spaces)  and sections of fibrations are assumed to be smooth.
Further, if an index exceeds its allowed range, the object labeled by this index is supposed to be zero.

\section{Parametrix possessing complexes}
 
To fix a terminology, we recall some notions from the theory of Hilbert $A$-modules.
Let  $A$ be a unital $C^*$-algebra. For a pre-Hilbert $A$-module $(U, (,)_U),$ let $| \mbox{ }|_U$ denote the associated  norm on $U$  defined by 
$|u| = \sqrt{|(u,u)_U|_A},$ $u\in U,$ where $| \mbox{ } |_A$ denotes the Banach norm on $A.$ If $(U, |\mbox{ }|_U)$ is a complete normed space,
 we call $(U,(,)_U)$ a Hilbert  $A$-module.
If $U$ is a 
pre-Hilbert $A$-module or a  Hilbert $A$-module, the $A$-valued map $(,)_U: U \times U \to A$ is called an $A$-product or a Hilbert $A$-product, respectively. 
For definiteness, we consider left pre-Hilbert $A$-modules, and the $A$-products are supposed to be conjugate linear in the first variable.
For pre-Hilbert $A$-modules $U,V,$ the space of continuous  $A$-module homomorphisms between $U$ and $V$ is denoted by $\mbox{Hom}_A(U,V)$ and its elements
 are called pre-Hilbert $A$-module homomorphisms. If $U, V$ are 
Hilbert $A$-modules, we omit the prefix pre. The notion of continuity is  meant with respect to the norms $| \mbox{ }  |_U$ and $| \mbox{ } |_V$. 
  Further, we say that $u,v\in U$ are orthogonal if $(u,v)_U=0.$
When we write a finite direct sum, the summands are supposed to be mutually orthogonal pre-Hilbert $A$-modules.
The adjoints  of maps between pre-Hilbert $A$-modules are always thought with 
respect to the considered $A$-products.
 Let us remark that there exist continuous maps  on a Hilbert $A$-module $(U, (,)_U)$ which are adjointable with respect to a
Hilbert space scalar product $U \times U \to \mathbb{C}$ on $U,$ the induced norm of which is  equivalent to the norm $| \mbox{ }|_U,$
 but  which do not posses an  adjoint with respect to  $(,)_U.$ For it, see, e.g., Solovyov, Troitsky \cite{TM}. A Hilbert $A$-module $U$ is called projective if  there exist $n \in \mathbb{N}_0$ and a Hilbert $A$-module $V$ such 
$A^n = U\oplus V,$ where $A^n$ denotes the direct sum of    $n$ copies of the standard  Hilbert $A$-module $A.$ It is called
finitely generated if there exists a finite number of elements $u_1, \ldots, u_k \in U$  so that 
for each $u\in U$ there exist elements $a_i \in A, i= 1,\ldots, k,$ such that $u=\sum_{i=1}^k a_i u_i.$  Let us notice,
 that the last mentioned Hilbert $A$-modules are sometimes called 
algebraically finitely generated.
For these notions, we refer the reader to Paschke \cite{Paschke} and Solovyov, Troitsky \cite{TM}.
 
Let us recall the following statement generalizing the rank-nullity theorem from linear algebra.

{\bf Theorem 1:} Let $U, V$ be  Hilbert $A$-modules and   $L \in \mbox{Hom}_A(U,V)$ be an adjointable map. If
the image of $L$ is closed, then also the image of $L^*$ is closed and $\mbox{Rng} \, L = (\mbox{Ker} \, L^*)^{\bot},$ $\mbox{Rng} \, L^* = (\mbox{Ker} \, L)^{\bot}$  
and $U=\mbox{Ker} \, L \oplus  \mbox{Rng} \, L^*.$

{\it Proof.} See Lance \cite{Lance} and a reference to the original proof there. $\Box$

{\bf Remark 1:} Under the same assumptions as in Theorem 1, the equality  $\mbox{Ker} \, {L}^* = (\mbox{Rng} \, {L})^{\bot}$ follows immediately.

{\bf Lemma 2}: Let $U, V, W$ be pre-Hilbert $A$-modules and 
$$\begin{xy} \xymatrix{
{U}   \ar[r]^D      & {V}   \ar[r]^{D'}  & {W}} 
\end{xy}$$ be a sequence of adjointable pre-Hilbert $A$-module homomorphisms.
Then for $\triangle = D'^* D' + DD^*,$ we have
$$\mbox{Ker} \, \triangle = \mbox{Ker} \, D'  \cap \mbox{Ker} \, D^*.$$ 

{\it Proof.} From the definition of $\triangle,$ we get the inclusion $\mbox{Ker} \, D' \cap \mbox{Ker} \, D^* \subseteq \mbox{Ker} \, \triangle.$ 
 For each $v\in V,$ we may write $(v, \triangle v)_V = (v, D'^*D' v+DD^*v)_V= (v,D'^*D'v)_V + (v, DD^*v)_V = (D'v,D'v)_W + (D^*v,D^*v)_U.$ 
 From this, the opposite inclusion follows using the positive definiteness of the $A$-products on $W$ and $U.$ $\Box$

To each complex $D^{\bullet}=(D_i, \Gamma^i)_{i\in \mathbb{N}_0}$ of pre-Hilbert $A$-modules and adjointable pre-Hilbert $A$-module homomorphisms, we attach the sequence 
$$\triangle_i = D_{i-1}D_{i-1}^*+D_{i}^*D_i, \mbox{  } i\in \mathbb{N}_0,$$ 
of {\it associated Laplacians}, where we assume $D_{-1}=0$  (according to the preamble).
 
{\bf Theorem 3:} Let $\Gamma^i,$ $i=1,\ldots, 5,$ be pre-Hilbert $A$-modules, and  let
$$\begin{xy} 
\xymatrix{
{\Gamma^1}   \ar[r]^{D_1}      & {\Gamma^2}   \ar[r]^{D_2}  & {\Gamma^3} \ar[r]^{D_3} & {\Gamma^4}  \ar[r]^{D_4} & {\Gamma^5} }
\end{xy}$$ be a complex of adjointable pre-Hilbert $A$-module homomorphisms. 
Suppose that for $i = 1, \ldots, 4,$ there exist elements  $g_i, p_i \in \mbox{Hom}_A (\Gamma^i,\Gamma^i)$  such that
\begin{eqnarray}
1_{|\Gamma^i} &=&   g_i  \triangle_i + p_i,\\
1_{|\Gamma^i} &=&   \triangle_i g_i + p_i \mbox{ and}\\
\triangle_i p_i &=& 0.
\end{eqnarray}
Then for $i = 1, 2, 3,$ 
\begin{eqnarray}
p_{i+1}D_i=0 \mbox{ and   } D_{i} g_{i} = g_{i+1} D_{i}.
\end{eqnarray}

{\it Proof.} 
We  use the relation $D_i p_i = 0$, $i=1,\ldots,4,$ repeatedly, which follows from Lemma 2 and from (3).
For $i = 1, 2, 3$ and $u \in \Gamma^i, $ we may write
\begin{eqnarray*}
p_{i+1}D_iu &=& p_{i+1}D_i (\triangle_i g_i u + p_i u)\\
           &=& p_{i+1}D_i \triangle_i g_i u + p_{i+1}(D_i p_i) u\\
           &=& p_{i+1}D_i (D_{i-1}D_{i-1}^* + D_i^* D_i) g_i u \\
           &=& p_{i+1}(D_i D_i^* )D_i g_i u \\
           &=& p_{i+1}\triangle_{i+1} D_i g_i u\\
           &=& (1-\triangle_{i+1}g_{i+1})\triangle_{i+1} D_i g_i u \\
           &=& \triangle_{i+1} D_i g_i u  - \triangle_{i+1}g_{i+1}\triangle_{i+1} D_i g_i u\\
           &=& \triangle_{i+1} D_i g_i u  - \triangle_{i+1}(1-p_{i+1}) D_i g_i u\\ 
           &=& \triangle_{i+1} D_i g_i u  - \triangle_{i+1} D_i g_i u + \triangle_{i+1}p_{i+1} D_i g_i u\\
           &=& 0,
\end{eqnarray*}
where we used $\triangle_{i+1}p_{i+1} =0,$ $i=1,2,3,$ in the last step. Notice that in the above computation (rows 2 and 5), we derived also  the relation $D_i\triangle_i = \triangle_{i+1}D_i,$ 
which we use in what follows. For $i = 1, 2,3$ and $u \in \Gamma^i,$ we have 
\begin{eqnarray*}
g_{i+1}D_i u &=& g_{i+1}D_i (\triangle_i g_i  + p_i) u \\
             &=& g_{i+1}D_i\triangle_i g_ i u + g_{i+1} (D_i p_i) u \\
             &=& g_{i+1}\triangle_{i+1}D_i g_i u \\
             &=& (1-p_{i+1})  D_i g_i u = D_ig_i u, 
\end{eqnarray*}
where at the last step, we used the relation $p_{i+1}D_i = 0$ derived above. $\Box$

{\bf Remark 2:}\begin{itemize} 
\item[1)] We call equations $(1)$ and $(2)$ the parametrix equations, and  any map $g_i$ satisfying  (1), (2) and (3) the pseudoinverse of 
$\triangle_i.$ In the special case $A=\mathbb{C},$ the name Green's functions is also used.
\item[2)] The second relation in row (4) says that  $g^{\bullet}=(g_i)_{i=1,2,3}$ is a chain map.
\end{itemize}

{\bf Remark 3:}\begin{itemize}
\item[1)] Notice that from $(1),$ we get  ${p_i}_{|\mbox{Ker} \, \triangle_i} = \mbox{1}_{|\mbox{Ker} \, \triangle_i},$ $i=1,\ldots,4.$ Using this fact and relation $(3),$ we see that
$p_i$ maps $\Gamma^i$ onto $\mbox{Ker} \,\triangle_i.$ 
\item[2)] Suppose now that  $(1),$ $(2)$ and $(3)$ are valid and the relation
\begin{eqnarray}
g_i p_i = 0
\end{eqnarray}
is satisfied for $i=1,\ldots, 4.$
Applying $p_i$ on $(2)$ from the right, we get  $p_i^2=p_i,$ i.e.,  $p_i$ is an idempotent in $\mbox{Hom}_A(\Gamma^i,\Gamma^i),$ $i=1,\ldots, 4.$
 

\end{itemize}

{\bf Definition 1:} 
We call a complex $D^{\bullet} = (D_i,\Gamma^i)_{i\in \mathbb{N}_0}$ of pre-Hilbert $A$-modules and  adjointable pre-Hilbert $A$-module homomorphisms  {\it parametrix possessing} if for each $i\in \mathbb{N}_0,$ there exist elements $g_i, p_i \in \mbox{Hom}_A(\Gamma^i, \Gamma^i)$ 
  such that relations (1), (2) and (3) are satisfied.

For a complex  $D^{\bullet}=(D_i,\Gamma^i)_{i\in \mathbb{N}_0}$ of pre-Hilbert $A$-modules and pre-Hilbert $A$-module homomorphisms, we consider the  cohomology groups
 $$H^i(D,A)= \frac{\mbox{Ker}(D_i:\Gamma^i \to \Gamma^{i+1})}{\mbox{Rng}(D_{i-1}: \Gamma^{i-1} \to \Gamma^i)}, \mbox{  } i \in \mathbb{N}_0,$$
and denote the nominator, the  cycles of 
the complex,  by $Z^i(D,A),$ and the space of boundaries, $\mbox{Rng} \, (D_{i-1}:\Gamma^{i-1} \to \Gamma^{i}),$   by $B^i(D,A).$
The cohomology groups come up with the canonical quotient $A$-module structure. 
We do not speak about a pre-Hilbert $A$-module structure on them because we do not   know whether $B^i(D,A)$ is a closed subspace of $\Gamma^i.$
Further, let us denote the kernel of $\triangle_i: \Gamma^i \to \Gamma^i,$ the $A$-module of {\it harmonic elements}, by $K^i(D,A).$

{\bf Theorem 4:} If $D^{\bullet}=(D_i, \Gamma^i)_{i\in \mathbb{N}_0}$ is a parametrix possessing complex, 
then for each $i\in \mathbb{N}_0,$  $$H^i(D,A) \simeq K^i(D,A) \mbox{ as $A$-modules}.$$

{\it Proof.}  Because $D^{\bullet}$ is parametrix possessing, there exist maps $p_i, g_i \in \mbox{Hom}_A(\Gamma^i, \Gamma^i)$ satisfying $(1),$ $(2)$ and $(3)$ for each
$i\in \mathbb{N}_0.$ 
 \begin{itemize}
\item[1)] Consider the  map
$\Phi_i:  K^i(D,A) \to H^i(D,A)$ given by $\Phi_ia = [a],$ $a\in K^i(D,A).$ This map is well defined. 
Indeed, due to Lemma 2, $a \in \mbox{Ker} \, \triangle_i$ implies $a \in \mbox{Ker} \, D_i.$
It is evident that $\Phi_i$ is an $A$-module homomorphism.

\item[2)] Let us consider the map $\Psi_i: H^i(D,A) \to K^i(D,A)$ given by
$\Psi_i[a] = p_ia$ for each $a \in Z^i(D,A).$
It follows from   Theorem 3 that $\Psi_i$ is well defined.
Indeed, for $a = D_{i-1}b$ we have $p_{i}a = p_i D_{i-1}b = 0$ due to the first relation in row (4).
Now, we prove that $\Psi_i$ is inverse to $\Phi_i.$ 
For $a \in K^i(D,A),$ we have $\Psi_i(\Phi_i a) = \Psi_i[a]=p_ia=a.$  
On the other hand, $\Phi_i(\Psi_i[b]) = \Phi_i(p_i b)=[p_i b ]$ for any $b\in Z^i(D,A).$
 The proof of the fact that $b$ and $\tilde{b}=p_i b$ are 
cohomologous proceeds  as follows. Using the parametrix equation (1), we get  $b - \tilde{b} = b - p_i b = (1-p_i) b = g_i\triangle_i b =
g_iD_{i-1}D^*_{i-1}b+g_iD_i^*D_ib = g_i D_{i-1} D^*_{i-1} b = D_{i-1}(g_{i-1}D_{i-1}^*b),$ 
where we used the  second equation in row (4) from Theorem 3 and the fact that $b \in \mbox{Ker} \, D_i.$
It is evident that $\Psi_i$ is an $A$-module homomorphism.
$\Box$
\end{itemize}

\section{$A$-elliptic pseudodifferential operators}

Let $p: E \to M^n$ be a locally trivial Banach bundle and $A$ be a unital $C^*$-algebra.
We call $p$ an {\it  $A$-Hilbert bundle} if 
\begin{itemize}
\item[1)] there exists a Hilbert $A$-module $(U, (,)_U)$ (the typical fiber of $p$),
\item[2)] for each $x\in M,$ the fiber $E_x= p^{-1}(x)$ is equipped with a Hilbert $A$-module structure, and it is isomorphic to $U$ as a Hilbert $A$-module,
\item[3)] the subset topology on $p^{-1}(x) \subseteq E$ is equivalent to the norm topology on $p^{-1}(x)$ given by $|\mbox{ }|_U,$ and
\item[4)] the local transition maps between the bundle charts are maps into the group $\mbox{Aut}_A(U)$ of Hilbert $A$-module automorphisms of $U.$
\end{itemize}

We call a smooth map $S: E_1 \to E_2$ between the total spaces of $A$-Hilbert bundles $p_1:E_1 \to M$ and $p_2: E_2 \to M$ with typical fibers $U_1$ and $U_2,$ respectively,
an $A$-Hilbert bundle homomorphism if it satisfies the equation  $p_2 \circ S = p_1$ (defining a bundle homomorphism),  and if it is a Hilbert $A$-module homomorphism in each fiber. 
Let $S$ be an $A$-Hilbert bundle homomorphism.
 We call an $A$-Hilbert bundle homomorphism $T:E_2 \to E_1$ adjoint to
$S$ if $T_{|p_2^{-1}(x)}$ is adjoint  to the Hilbert $A$-module homomorphism $S_{|p_1^{-1}(x)} : p_1^{-1}(x) \to p_2^{-1}(x)$ for each $x\in M.$ 
In this case, we write $T=S^*.$ 
An $A$-Hilbert bundle is called  {\it finitely generated projective} $A$-Hilbert bundle if its typical fiber is a finitely generated projective Hilbert $A$-module.

The vector space $\Gamma(E)$ of smooth sections of $p$ carries a  structure of an $A$-module given by $(a s)(x) = a(s(x)),$ where $a\in A,$  $s\in \Gamma(E)$ and  $x\in M.$ 
For a compact manifold $M$  and a Riemannian metric $g$ on $M,$ we fix  an associated volume element on $M$  and denote it by $|\mbox{vol}_g|.$ 
The volume element
 induces an  $A$-product on $\Gamma(E)$  by the formula
$$(s,s')_{\Gamma(E)} = \int_{x \in M} (s(x),s'(x))_U |\mbox{vol}_g(x)| \in A,$$ where 
$s, s' \in \Gamma(E).$ At the right hand side, we consider  the Bochner integral. 
 In this way, $\Gamma(E)$   gains a structure of a pre-Hilbert $A$-module. We denote the induced norm by $| \mbox{ }|_{\Gamma(E)}.$
 For each $t\in \mathbb{Z},$ one defines a further $A$-product $(,)_t$ on $\Gamma(E)$ by setting
$$(s,s')_t = \int_{x\in M} ((1+\triangle_g)^t s(x),s'(x))_U |\mbox{vol}_g(x)| \in A,$$ where
$\triangle_g$ is  minus the Laplace-Beltrami operator of $(M, g)$ and $s, s' \in \Gamma(E).$ The induced norm
will be denoted by $| \mbox{ }|_t,$ and the completion of $\Gamma(E)$ with respect to $| \mbox{ }|_t$ by 
$W^t(E).$ In particular, $\Gamma(E)$ is dense in $W^t(E)$.
Due to the construction,  $W^t(E)$ together with the extended $A$-product $(,)_t$ form a Hilbert $A$-module. We call any
$(W^t(E), (,)_t)$ a {\it Sobolev type completion} of $\Gamma(E).$  
See Fomenko, Mishchenko \cite{FM} for more on the introduced $A$-products.
 
Because ${| \mbox{ } |_0}_{|\Gamma(E)} = | \mbox{ } |_{\Gamma(E)},$ the Hilbert $A$-module $W^0(E)$ coincides
 with the completion of $\Gamma(E)$ with respect to $| \mbox{ } |_{\Gamma(E)}.$
For $t\in \mathbb{Z}$ and $k\geq t>0,$ it is not difficult to see that
 $$\Gamma(E) \subseteq \mathcal{C}^k(E) \subseteq W^{t}(E) \subseteq W^{t-1}(E).$$
A section of $E$ belongs to $\mathcal{C}^k(E)$ if and only if it is at least $k$ times continuously differentiable.
To be more precise, let us notice that we consider the elements of $\Gamma(E), C^k(E)$ and $W^t(E)$ modulo the relation of being zero almost everywhere.
For any $O$ open in $M,$ we  use the symbol $\mathcal{C}^k(O,E)$ to denote the space of the appropriate  sections of the restricted bundle $p_O: p^{-1}(O) \to O.$

 Let us briefly recall a definition of the Fourier 
transform of  local sections of  $A$-Hilbert bundles over compact Riemannian manifolds. 
For a local chart $(O \subseteq \mathbb{R}^n, \phi: O \to \phi(O) \subseteq M^n),$ the Fourier transform of a section 
$s$ of $p_{\phi(O)},$ which has its support in a compact subset of $\phi(O),$
  is  defined by $$\hat{s}(q)=  \int_{x\in \phi(O)} e^{-2\pi \imath \langle q, x \rangle }s(x) \mbox{d}x,$$
where   $q \in \mathbb{R}^n$  and $\langle, \rangle$ denotes the appropriate Euclidean product on $\mathbb{R}^n$. We use the same notation for 
a section as for its coordinate expression, hoping this causes no confusion.
  Globally, one has to choose an atlas and a subordinate 
partition of unity, make the local Fourier transforms, and apply the appropriate gluing process at the end.
Let us notice that the Fourier transform depends on the choice of a particular partition, but it exists if the underlying Riemannian manifold is 
compact independently of this choice.
As in the classical case, one can show that the norm $| \mbox{ }|_t$ is equivalent to the norm 
$||s||_t=[\int_{q \in \mathbb{R}^n}|\hat{s}(q)|^2_{\Gamma(E)}(1 + |q|^2)^t \mbox{d}q]^{\frac{1}{2}} \in \mathbb{R}_0^+,$ where $s\in \Gamma(E),$
 $|q|=\langle q,q \rangle^{\frac{1}{2}},$  and  $q\in \mathbb{R}^n.$ 
To obtain a global version of this formula, one shall apply the same procedure as in the case of the Fourier transform.

 In the following lemma, an $A$-Hilbert bundle analogue of the Sobolev embedding theorem is proved. 

{\bf Lemma 5:} Let $p: E\to M^n$ be an  $A$-Hilbert bundle 
over a compact manifold $M.$  Then  for  each $t > \lfloor \frac{n}{2} \rfloor +1 $ and  $ 0\leq  k < t - \lfloor \frac{n}{2}\rfloor - 1,$ 
$$s\in W^t(E) \mbox{ implies  } s \in \mathcal{C}^k(E).$$

{\it Proof.} Let $t > \lfloor \frac{n}{2} \rfloor + 1$ and $0\leq  k < t - \lfloor \frac{n}{2}\rfloor - 1.$
 If $s$ is a section of $p_O,$ $O\subseteq M,$ we denote its coordinate expression by $s$ as well.
  For $\alpha \in \mathbb{N}_0^n$  and $|\alpha| \leq k,$ we show that
 $\partial^{\alpha}s \in \mathcal{C}^0(O,E).$ 
Obviously, it is sufficient to prove that
$\int_{q \in \mathbb{R}^n} \hat{s}(q) q^{\alpha} e^{2 \pi \imath \langle x,q \rangle} \mbox{d}q$ 
converges for all $x\in O.$ Because $s \in W^t(O,E),$ we know that $||\hat{s}||_t = [\int_{q\in \mathbb{R}^n}
 |\hat{s}(q)|_{\Gamma(E)}^2 $ $(1+|q|^2)^t \mbox{d}q]^{\frac{1}{2}} < \infty.$
Let us compute
\begin{eqnarray*}
&&\int_{q \in \mathbb{R}^n} |\hat{s}(q)|_{\Gamma(E)} |q^{\alpha}| \mbox{d}q \leq\int_{q \in \mathbb{R}^n} |\hat{s}(q)|_{\Gamma(E)}  |q|^{|\alpha|}  \mbox{d}q =\\
&=& \int_{q \in \mathbb{R}^n} |\hat{s}(q)|_{\Gamma(E)}(1+|q|^2)^{t/2}\frac{|q|^{|\alpha|}}{(1+|q|^2)^{t/2}} \mbox{d}q \leq\\
&\leq& \left(\int_{q \in \mathbb{R}^n} |\hat{s}(q)|^2_{\Gamma(E)}(1+|q|^2)^t \mbox{d} q \right)^{1/2}   
\left(\int_{q \in \mathbb{R}^n} \frac{|q|^{2|\alpha|}}{(1+|q|^2)^t} \mbox{d} q\right)^{1/2}=\\
&=&||\hat{s}||_t  \left(\int_{q \in \mathbb{R}^n} \frac{|q|^{2|\alpha|}}{(1+|q|^2)^t} \mbox{d}q \right)^{1/2},
\end{eqnarray*}
where we used the Cauchy-Schwartz inequality in $L^2(\mathbb{R}^n)$ at the second last step.
 Using   polar coordinates in $\mathbb{R}^n$, we see that the  finiteness of the last written integral  is equivalent to the finiteness
of $\int_{0}^{+\infty} r^{2|\alpha|+n-1}(1+r^2)^{-t} \mbox{d}r.$ Near $r=0,$ the integrand 
is a bounded continuous function. At the infinity, the integrand behaves like
$r^{2|\alpha| + n - 2t}$ which is  integrable over $(C, +\infty)$ for each  $C>0$ if and only if $2 |\alpha| + n - 2 t  < -1.$ 
Thus, for each $|\alpha| < t - \frac{n}{2} -\frac{1}{2},$ the integral $\int_{q \in \mathbb{R}^n} \hat{s}(q) q^{\alpha} e^{2 \pi \imath \langle x,q\rangle} \mbox{d}q$  converges due to the absolute convergence of 
the Bochner integral.    Therefore $\partial^{\alpha}s \in \mathcal{C}^0(O,E)$ 
for each $|\alpha|\leq k < t -\lfloor\frac{n}{2} \rfloor -1.$ $\Box$

A  reference for   proofs of the statements in this paragraph is Solovyov, Troitsky \cite{TM}.
For $A$-Hilbert bundles $p_1:E_1\to M^n$  and $p_2:E_2 \to M,$  one defines a symbol of a pseudodifferential operator of order $r \in \mathbb{Z}.$
Besided other properties, the symbol is an adjointable $A$-Hilbert bundle homomorphism.
Using the Fourier transform, to each symbol of a pseudodifferential operator of order $r$,   one  associates a pseudodifferential operator of order $r$ and vice versa. 
Let us notice that any  pseudodifferential operator $D: \Gamma(E_1) \to \Gamma(E_2)$ of order $r$ is, in  particular, an adjointable pre-Hilbert $A$-module
 homomorphism  between the   pre-Hilbert $A$-modules $\Gamma(E_1)$ and $\Gamma(E_2)$. The order of $D^*$ equals that of the operator  $D.$
Further, it is  known   that an arbitrary $t\in \mathbb{Z},$ each pseudodifferential operator $D: \Gamma(E_1) \to \Gamma(E_2)$  of order $r$ admits  a unique Hilbert $A$-module
homomorphism $D^t: W^t(E_1) \to W^{t-r}(E_2)$ extending  $D$ and that  any of these extensions is adjointable. Further, the symbol of $D^*$ equals  (up to a constant non-zero multiple) to the adjoint of the symbol of $D.$
 
 One calls a pseudodifferential operator $K:\Gamma(E_1) \to \Gamma(E_2)$  {\it  $A$-elliptic} if  
for each $x \in M$ and $\xi \in T_x^*M \setminus \{0\},$
the symbol $$\sigma(x,\xi):(E_1)_x \to (E_2)_x $$ of $K$ at $(x,\xi)$ is a Hilbert $A$-module isomorphism.
Further, let $U, V$ be Hilbert $A$-modules, and the decompositions $U = U_0 \oplus U_1$ and $V = V_0\oplus V_1$ hold, where $U_1$ and $V_1$ are Hilbert $A$-modules,
 and $U_0$ and $V_0$ are finitely generated projective 
Hilbert $A$-modules. Whenever $K_i \in \mbox{Hom}_A(U_i,V_i),$ $i=0,1,$ and $K_1$ is an isomorphisms, then
$K = K_0 \oplus K_1$ is called an $A$-Fredholm operator.
 It is immediately   seen that $\mbox{Ker} \, K= \mbox{Ker} \, K_0$  and $\mbox{Coker}\, K = V_0/\mbox{Rng}\, K_0.$
Let us notice that in general, the image of an $A$-Fredholm operator  is not closed in contrary to the case of $A=\mathbb{C}$.
 
Now, we recall the theorem of Fomenko and Mishchenko on the existence of smoothing pseudoinverses  mentioned in the Introduction.

{\bf Theorem 6:} Let $p: E \to M$  be a finitely 
generated projective $A$-Hilbert bundle over a compact manifold $M,$ and  $K: \Gamma(E) \to \Gamma(E)$ be an $A$-elliptic   operator of order $r$. 
Then for each $t \in \mathbb{Z},$ the extension $K^t$ of $K$ is an  $A$-Fredholm operator and there exists a Hilbert $A$-module homomorphism $g^{t-r}: W^{t-r}(E)\to W^t(E)$
satisfying $g^{t-r}K^t - 1 : W^t(E) \to W^{t+1}(E).$
 
{\it Proof.}
See Solovyov, Troitsky \cite{TM} (Theorems 2.1.142 and 2.1.146). $\Box$
 
Now, it is quite easy to prove the regularity of  $A$-elliptic operators.

{\bf Theorem 7:}
Let $p: E \to M$ be a finitely generated projective $A$-Hilbert bundle over a compact manifold $M$  and $K: \Gamma(E) \to \Gamma(E)$
be an $A$-elliptic   operator of order $r$. Then  for each $t \in \mathbb{Z},$  the equality
$\mbox{Ker} \, K^t = \mbox{Ker} \, K \subseteq \Gamma(E)$ holds.

{\it Proof.}  The inclusion $\mbox{Ker}\, K \subseteq \mbox{Ker} \, K^t$ is obvious. For $t \in \mathbb{Z},$
let us prove that $\mbox{Ker} \, K^t \subseteq \mbox{Ker}\, K.$
Due to Theorem 6, there exists   a map $g^{t-r} : W^{t-r}(E) \to W^{t}(E)$
such that $g^{t-r} K^t - 1 : W^{t}(E) \to W^{t+1}(E).$
For $s \in W^{t}(E),$ we may write 
\begin{eqnarray*}
s  &=&  (g^{t-r} K^t)s   -  (g^{t-r}K^t)s + s  \\
   &=&   g^{t-r} (K^t s) -  (g^{t-r}K^t - 1)s.
\end{eqnarray*}
Assuming $K^t s = 0$ for $s\in W^t(E),$ we get $s=g^{t-r} (K^t s)- (g^{t-r} K^t - 1)s = - (g^{t-r} K^t -1) s \in W^{t+1}(E)$ (Theorem 6).
By induction, $s \in W^t(E)$ for each $t\in \mathbb{Z}.$  Using Lemma 5, we obtain $s\in \bigcap_{k\in \mathbb{N}_0} \mathcal{C}^k(E) = \Gamma(E).$ $\Box$

Using the previous theorem, we prove a smooth version of Theorem 6.

{\bf Theorem 8:} Let $p: E \to M$ be a finitely generated projective 
$A$-Hilbert bundle over a compact manifold $M,$ and  $K: \Gamma(E) \to \Gamma(E)$ be a self-adjoint $A$-elliptic   operator of order $r$. 
If $\mbox{Rng} \, K^r$ is closed in $W^{0}(E),$ then there exist
pre-Hilbert $A$-module homomorphisms  $P: \Gamma(E) \to \Gamma(E)$ and $G: \Gamma(E) \to \Gamma(E)$ such that
 $$1_{|\Gamma(E)} = GK + P, \mbox{  } 1_{|\Gamma(E)} =  K G + P \mbox{ and } K P = 0.$$
{\it Proof.}
Notice that $K^r:W^r(E) \to W^{0}(E)$ and ${K^r}^*: W^{0}(E) \to W^r(E) \subseteq W^{-r}(E).$
\begin{itemize}
\item[1)]For each $a, b \in \Gamma(E),$ we have $((K^r)^*a ,b)_r = (a,K^r b)_0 = (a,Kb)_{\Gamma(E)} = (a,K^*b)_{\Gamma(E)} = (Ka, b)_{\Gamma(E)}.$
Summing up,  $((K^r)^*a,b)_r = (Ka,b)_{\Gamma(E)}$ for any $a,b \in \Gamma(E).$
For each  $a\in \mbox{Ker}\, K \subseteq \Gamma(E)$ (Theorem 7), we thus get $((K^r)^*a,b)_r=0$ for any $b\in \Gamma(E).$ 
Because $(, )_r$ is continuous and $\Gamma(E)$ is dense in $W^r(E),$ we obtain  $((K^r)^*a,b)_r=0$ for each $b\in W^r(E).$
Thus from the non-degeneracy of $(,)_r,$  ${K^r}^*a=0,$  and consequently, $\mbox{Ker} \, K \subseteq \mbox{Ker}\, {K^r}^*.$  
If $a \in \mbox{Ker} \, {K^r}^* \cap \Gamma(E),$ then from the previous computation and the non-degeneracy of $(,)_0$, we get $a\in \mbox{Ker} \, K.$
Therefore  $\mbox{Ker} \, K = \mbox{Ker}\, {K^r}^* \cap \Gamma(E).$ Because ${K^r}^*$ is also an $A$-elliptic operator, Theorem 7 applies, and
 we may omit the intersection with
$\Gamma(E)$ in the previous expression, obtaining
\begin{eqnarray}
\mbox{Ker} \, K = \mbox{Ker}\, {K^r}^*.
\end{eqnarray}

%
%
%
%

\item[2)]Because the image of $K^r: W^r(E) \to W^{0}(E)$ is   closed and $K^r$ is adjointable, we have
$W^r(E)=\mbox{Ker}\, K^r \oplus\mbox{Rng} \, {K^r}^*$  due to Theorem 1. Therefore, the projection  $p_{\mbox{Ker} \, K^r}: W^r(E) \to
\mbox{Ker} \,K^r$ from $W^r(E)$ onto $\mbox{Ker} \, K^r$ is a well defined Hilbert $A$-module homomorphism.
Because $\mbox{Rng} \, K^r  = (\mbox{Ker} \, {K^r}^*)^{\bot} \subseteq W^{0}(E)$ (Theorem 1), there exists a bijective Hilbert $A$-module homomorphism
$$\delta  = K^r_{|(\mbox{Ker} \, K^r)^{\bot}}: (\mbox{Ker} \, K^r)^{\bot} \to (\mbox{Ker} \, {K^r}^*)^{\bot}.$$ 
Due to  Banach's open mapping theorem, 
 its inverse
$$\gamma:(\mbox{Ker} \, {K^r}^*)^{\bot} \to (\mbox{Ker}\,  K^r)^{\bot}$$ is continuous.
 As any inverse of an $A$-module homomorphism, the map $\gamma$ is an $A$-module homomorphism as well.
Extending $\gamma$ by zero on $\mbox{Ker} \, {K^r}^*=(\mbox{Rng}\, K^r)^{\bot}$ (Remark 1),  we get a Hilbert 
$A$-module homomorphism $$\tilde{\gamma}: W^0(E) \to (\mbox{Ker} \, K^r)^{\bot}
 \subseteq W^{r}(E).$$ 

\item[3)] Due to the construction, we obtain a parametrix type equation $1_{|W^r(E)} = \tilde{\gamma} K^r + p_{\mbox{Ker} \, K^r}$ for
$K^r.$  
Indeed, for $a\in \mbox{Ker}\, K^r,$ we have $\tilde{\gamma} K^r a + p_{\mbox{Ker} \, K^r} a = 0 + a = a.$ 
For $a\in (\mbox{Ker}\, K^r)^{\bot},$ we get $\tilde{\gamma} K^r a + p_{\mbox{Ker} \, K^r} a =  \gamma K^r a + 0 = a.$

Let us denote the restriction to $\Gamma(E)$ of 
$p_{\mbox{Ker} \, K^r}$ by $P.$ Thus $P: \Gamma(E) \to \mbox{Ker} \, K^r.$ Because $\mbox{Ker} \, K^r = \mbox{Ker} \, K \subseteq \Gamma(E)$ (Theorem 7),
we have a pre-Hilbert $A$-module homomorphism $P:\Gamma(E) \to \Gamma(E)$ at our disposal.

Further, let us set $G=\tilde{\gamma}_{|\Gamma(E)}.$ 
For $a \in \mbox{Ker} \, {K^r}^* \subseteq  \Gamma(E),$ we have $G a= \tilde{\gamma} a = 0.$
For $a \in (\mbox{Ker} \, {K^r}^*)^{\bot} \cap \Gamma(E),$ there exists $b \in (\mbox{Ker} \, K^r)^{\bot}$ (item 2) such that
$a = \delta b.$ Since $a \in \Gamma(E)$ and $\delta$ is a restriction of an extension of a pseudodifferential operator of finite order,   $b\in \Gamma(E).$
We have $Ga = \gamma a = \gamma \delta b = b.$ Especially, $Ga \in \Gamma(E).$ Thus, 
$G: \Gamma(E) \to \Gamma(E).$
Now, we may restrict the parametrix type equation  for $K^r$ to the space $\Gamma(E)$ to obtain   $1_{|\Gamma(E)} = G K + P.$

\item[4)] Similarly as above, we prove that the parametrix type equation $1_{|W^{0}(E)}=K^r\tilde{\gamma} + p_{\mbox{Ker} \, {K^r}^*}$ holds,
 where $p_{\mbox{Ker} \, {K^r}^*}$ denotes the
 projection from $W^0(E)$ onto
$\mbox{Ker} \, {K^r}^*.$
Indeed, for $a\in \mbox{Ker} \, {K^r}^*,$ we obtain $K^r \tilde{\gamma} a + p_{\mbox{Ker} \, {K^r}^*} a = 0 + a =a.$
For $a\in (\mbox{Ker}\, {K^r}^*)^{\bot},$ we have $K^r \tilde{\gamma} a + p_{\mbox{Ker} \, {K^r}^*} a = a + 0 = a.$
Using formula $(6)$ from item 1, we see that $(p_{\mbox{Ker} {K^r}^*})_{|\Gamma(E)} = (p_{\mbox{Ker} \, K^r})_{|\Gamma(E)}.$ Thus restricting 
$1_{|W^{0}(E)}=K^r\tilde{\gamma} + p_{\mbox{Ker} \, {K^r}^*}$ to $\Gamma(E),$ we get $1_{|\Gamma(E)}=K G + P.$

\item[5)] The equation $K P = 0$ follows from  the definition of $P$ and the inclusion 
$\mbox{Ker} \, K^r \subseteq \Gamma(E)$ (Thm. 7).
\end{itemize}
$\Box$

{\bf Remark 4:} Under the conditions of Theorem 8 and from the constructions in its proof, we get that  $GP=0$
 and $P^2=P.$ See also Remark 3 item 2.

\section{Complexes of $A$-elliptic pseudodifferential operators}

In this section, we  focus our attention at complexes of  pseudodifferential operators rather than at one pseudodifferential operator only. 
Let $M$ be a manifold, $E^{\bullet}=(p_i: E^i \to M)_{i \in \mathbb{N}_0}$ be 
a sequence of $A$-Hilbert bundles over $M,$
and $D^{\bullet}=(D_i, \Gamma(E^i))_{i \in \mathbb{N}_0}$ be a complex of  pseudodifferential 
operators acting between   sections of the appropriate bundles, $A$ being a unital $C^*$-algebra.
Let $\sigma_i$ be the symbol of the pseudodifferential operator $D_i: \Gamma(E^i) \to \Gamma(E^{i+1}),$ $i \in \mathbb{N}_0.$ 
For each $i \in \mathbb{N}_0,$ we define  $\sigma_i'$ to be the restriction of the $A$-Hilbert bundle homomorphism 
$\sigma_i$ to $\pi^*(E^i),$ where $\pi: T^*M\setminus 0 \to M$ is the foot point projection from  the cotangent bundle with the image of the zero section removed onto the manifold $M$.
Let us denote the resulting restricted complex consisting of $A$-Hilbert bundle homomorphisms $\sigma_i',$ $i\in \mathbb{N}_0,$ by ${\sigma^{\bullet}}'.$
 

{\bf Definition 2:}
A complex $D^{\bullet}=(D_i,\Gamma(E^i))_{i\in \mathbb{N}_0}$ of pseudodifferential operators is called an {\it $A$-elliptic}  complex if 
  ${\sigma^{\bullet}}'$ is an exact complex in the category of $A$-Hilbert bundles.
 
The following lemma is a  generalization to the case of Hilbert $A$-modules of a  reasoning used in the classical Hodge theory.
For the classical case, the reader is referred, e.g., to Wells \cite{Wells}.
 
{\bf Lemma 9:} Let $U, V, W$ be Hilbert $A$-modules and  $\sigma$ and ${\sigma'}$ be adjointable  Hilbert
$A$-module homomorphisms. If
$$\begin{xy} \xymatrix{
{U}   \ar[r]^\sigma    & {V}   \ar[r]^{\sigma'}  & {W}} 
\end{xy}$$
is an exact sequence and the image of $\sigma'$ is closed, then $\Sigma = \sigma \sigma^*  + {\sigma'}^*{\sigma'}$ is a Hilbert $A$-module automorphism of $V$.

{\it Proof.} 
 Due to the assumption, the Hilbert $A$-module homomorphisms  $\sigma^*$ and ${\sigma'}^*$ exist.
Therefore, $\Sigma = \sigma \sigma^* + {\sigma'}^*{\sigma'}$ is a well defined Hilbert $A$-module endomorphism. It is immediately seen that $\mbox{Rng}\,{\sigma'}^* \subseteq \mbox{Ker}\,\sigma^*.$
 Indeed, for any $v \in \mbox{Rng}\, {\sigma'}^*$  there exists an element $w\in W$ such that
$\sigma'^*w=v.$ Thus $\sigma^* v = (\sigma^*{\sigma'}^*) w= (\sigma'\sigma)^*w=0.$ 

Now, we prove the injectivity of $\sigma \sigma^*$  restricted to $\mbox{Rng} \, \sigma.$ Suppose there exists an element $u \in U$  such that
$v=\sigma u$ satisfies  $\sigma \sigma^* v =0.$ We may write $0 = (\sigma \sigma^* v, v)_V = (\sigma^*v, \sigma^*v)_U$  which  implies $\sigma^*v =0.$ Thus
 $0 = (\sigma^*v, u)_U = (v, \sigma u)_V = (v, v)_V$
which implies $v=0.$   Now, we prove that ${\sigma'}^* {\sigma'}$ 
is injective on $\mbox{Rng} \, {\sigma'}^*.$
Indeed, let $v = {\sigma'}^* w$ for  an element  $w\in W$ and ${\sigma'}^* {\sigma'} v =0.$ Then ${\sigma'}^*{\sigma'} v = 0$ 
implies that $(\sigma' v, \sigma' v)_W = (v, {\sigma'}^* {\sigma'} v)_V = 0,$ which  in turn implies  $\sigma' v = 0.$ 
 We may  compute $0 = ({\sigma'} v, w)_W = ( v, {\sigma'}^* w)_V =  (v,v)_V$  from  which
$v = 0$ follows. Therefore ${\sigma'}^*{\sigma'}$ is injective on $\mbox{Rng}\, {\sigma'}^*.$

Further, due to Theorem 1,
we have $V = \mbox{Ker}\, {\sigma'} \oplus \mbox{Rng} \, {\sigma'}^*.$ Because 
$\mbox{Rng} \, \sigma = \mbox{Ker} \, {\sigma'}$, we may write 
$V = \mbox{Rng} \, \sigma \oplus \mbox{Rng} \,  {\sigma'}^*.$ Since $\mbox{Rng} \, \sigma \subseteq \mbox{Ker}\,  \sigma'$
and $\mbox{Rng} \, {\sigma'}^* \subseteq \mbox{Ker}\, \sigma^*,$ we get
$\Sigma_{|\mbox{Rng} \, \sigma} = \sigma \sigma^*_{|\mbox{Rng} \, \sigma}$ 
and $\Sigma_{|\mbox{Rng} \, {\sigma'}^*} = {\sigma'}^* \sigma'_{|\mbox{Rng} \, {\sigma'}^*}.$ Suppose there exists  $v = (v_1,v_2) \in
 \mbox{Rng}\, {\sigma} \oplus \mbox{Rng} \, {\sigma'}^* $ such that $\Sigma v =0.$ Then $0 = \Sigma v = \sigma\sigma^* v_1 + {\sigma'}^* \sigma' v_2.$
Because $\mbox{Rng} \, \sigma \cap \mbox{Rng} \, {\sigma'}^* = 0,$ we have $\sigma\sigma^* v_1 = {\sigma'}^*\sigma' v_2 = 0.$ Due to the
 injectivity of the appropriate restrictions of the maps $\sigma\sigma^*$ and ${\sigma'}^* \sigma'$, 
the last two equations imply
$v_1 = v_2 = 0.$ Thus $\Sigma$  is injective on $V.$ 

Since the images of $\sigma$ and $\sigma'$ are closed, the images of $\sigma\sigma^*$ and ${\sigma'}^*\sigma'$ are closed as well. 
(For it, see the proof of Theorem 3.2 in Lance \cite{Lance}.)
Because $\sigma\sigma^*$ and ${\sigma'}^* {\sigma}'$ have closed images, their sum $\Sigma$ has closed image as well. Obviously, $\Sigma$ is adjointable.
Due to Theorem 1, $\Sigma$ is surjective because
$V=\mbox{Rng}\,(\Sigma^*)\oplus \mbox{Ker}\,(\Sigma) = \mbox{Rng}\,(\Sigma^*) = \mbox{Rng}\,(\Sigma).$ 
$\Box$

{\bf Remark 5:} Going through the proof of the previous theorem, we  see   that under the same conditions on $\sigma$ and $\sigma',$ we obtain that 
the map $\lambda \sigma^*\sigma +  \mu \sigma\sigma^*$ is a  Hilbert $A$-module automorphism for all $\lambda, \mu \in \mathbb{C} \setminus \{0\}$ as well.

{\bf Corollary 10:} Let $M$  be manifold and $D^{\bullet}=(D_i, \Gamma(E^i))_{i\in \mathbb{N}_0}$ be an $A$-elliptic  complex 
in $A$-Hilbert bundles over $M$. 
Then for each $i\in \mathbb{N}_0,$ the associated Laplacian $\triangle_i$ is  $A$-elliptic.

{\it Proof.}  
Because $D^{\bullet}$ is an $A$-elliptic  
complex, the restricted symbol complex ${\sigma^{\bullet}}'$ 
is exact.
 From the exactness of ${\sigma^{\bullet}}',$ one concludes 
that the images of the Hilbert $A$-module homomorphisms $\sigma_{i}'$  evaluated at any $(x,\xi) \in T^*M \setminus 0$ 
are closed. Further,   they are adjointable, being the symbols
of pseudodifferential operators. 
Using Lemma 9,  the maps $\Sigma_i=\sigma_i^*\sigma_i + \sigma_{i-1}\sigma_{i-1}^*,$ 
$i\in \mathbb{N}_0,$  are Hilbert $A$-module automorphisms when evaluated at  $(x,\xi) \in T^*M \setminus 0.$ 
Consequently, the associated Laplacians $\triangle_i,$ $i \in \mathbb{N}_0,$ of $D^{\bullet}$ are $A$-elliptic operators.  
$\Box$


We denote the order of the associated Laplacian $\triangle_i$  by $r_i,$ $i \in \mathbb{N}_0.$
Further, we fix the topologies on the spaces which appear in the definition of the cohomology groups of a complex $D^{\bullet} = (D_i, 
\Gamma(E^i))_{i\in \mathbb{N}_0}$ of pseudodifferential operators.
We consider the spaces $\Gamma(E^i)$   with the topology given by the norm $| \mbox{ } |_{\Gamma(E^i)}$
and the spaces $Z^i(D,A), B^i(D,A) \subseteq \Gamma(E^i)$   with the    subset topologies. Notice that $\Gamma(E^i)$ is not a Banach space  in general.
We assume each cohomology group $H^i(D,A)$ to be equipped with the quotient topology. 
 
 
{\bf Theorem 11:} Let $D^{\bullet} = (D_i, \Gamma(E^i))_{i\in \mathbb{N}_0}$ be an $A$-elliptic  complex in finitely 
generated projective $A$-Hilbert bundles over  a compact manifold $M$.   
If the image of $\triangle_i^{r_i}$ is closed in $W^{0}(E^i),$ $i\in \mathbb{N}_0$, then for each $i \in \mathbb{N}_0$ 
\begin{itemize}
\item[1)]  the cohomology group $H^i(D,A)$ of $D^{\bullet}$ 
is a finitely generated $A$-module,
\item[2)] the image of $D_i$ is closed in the norm topology on $\Gamma(E^i),$ and
\item[3)] $H^i(D,A)$ is a Banach topological vector  space.
\end{itemize}

{\it Proof.} Let us fix a non-negative integer $i \in \mathbb{N}_0.$
According to Corollary 10, the Laplacian $\triangle_i$ is an $A$-elliptic   operator. Since it is self-adjoint  and 
$\mbox{Rng} \, \triangle_i^{r_i} \subseteq W^{0}(E_i)$ is closed,  $D^{\bullet}$ is a parametrix possessing complex according to Theorem 8.
Let us denote   the projection of $\Gamma(E^i)$ onto $K^i(D,A)$ from  Theorem 8 ($E=E^i,$ $K =\triangle_i$) by $P_i.$

\begin{itemize}
\item[1)] 
 Because $D^{\bullet}$ is a parametrix possessing complex, $K^i(D,A)$ is $A$-linearly isomorphic to $H^i(D,A)$ (Theorem 4).
Since $\triangle_i$ is $A$-elliptic, $\mbox{Ker} \, \triangle_i^{r_i} = K^i(D,A)$ due to the Theorem 7.
From Theorem 6, we know that $\triangle_i^{r_i}$ is $A$-Fredholm and thus, especially, $\mbox{Ker} \, \triangle^{r_i}_i$ is finitely generated as an $A$-module. Therefore 
$H^i(D,A)$ is finitely generated over $A$ as well.

\item[2)] Now, we prove that $H^i(D,A)$ is a Hausdorff space.
It is  obvious that the  maps
$$\Phi_i:K^i(D,A) \to H^i(D,A) \mbox{  and  } \Psi_i:H^i(D,A) \to K^i(D,A)$$ defined  by $\Phi_i a=[a]$ and
 $\Psi_i[a] =P_ia,$ $a\in Z^i(D,A)$   are continuous. For the fact that they are mutually inverse, see item 2 in the proof of Theorem 4.
Thus, $K^i(D,A)$ is  homeomorphic to $H^i(D,A).$  
Because   $K^i(D,A)$ is a  Hausdorff space, being a subspace of the Hausdorff space $\Gamma(E^i),$   $H^i(D,A)$ is a Hausdorff topological space as well. Since $H^i(D,A) = Z^i(D,A)/B^i(D,A)$ and $Z^i(D,A)$ is closed in $\Gamma(E^i)$, being the kernel of a continuous 
map, the space $B^i(D,A)$ of boundaries in $\Gamma(E^i)$ is closed in $\Gamma(E^i)$ as well. 
 
\item[3)] Proving that $K^i(D,A)$ is a Banach space with respect to the norm topology on $\Gamma(E^i)$ we are done due to item 2 of this proof.  
Let us consider the subset topology on $\mbox{Ker} \, \triangle_i^0$ inherited from the topology
of  $(W^0(E^i), |\mbox{ } |_0).$  
Because $|\mbox{ }|_{\Gamma(E^i)} = {| \mbox{ }|_0}_{|\Gamma(E)},$ any sequence in $K^i(D,A)$ which is cauchy with respect to  $| \mbox{ } |_{\Gamma(E^i)}$
 is also cauchy with respect to $| \mbox{  } |_0$  and thus, it has a limit in $W^{0}(E^i)$ since  $W^0(E^i)$ is complete.
 Because $K^i(D,A)$ equals $\mbox{Ker} \, (\triangle_i^{0}: W^0(E^i) \to W^{-r_i}(E^i))$ (Thm. 7), which is closed in $W^{0}(E^i)$
being the kernel of a continuous map,
the mentioned limit belongs also to $K^i(D,A).$ Since the norms $| \mbox{ } |_{\Gamma(E^i)}$  and $| \mbox{  } |_0$ (the latter restricted to $\Gamma(E^i))$ coincide,
 the limit is also a limit
with respect to the norm $| \mbox{  }|_{\Gamma(E^i)}.$  
 \end{itemize}
$\Box$


\end{document}